\documentclass{article}
\usepackage{epsfig}

\newtheorem{theorem}{Theorem}[section]
\newtheorem{definition}{Definition}

\begin{document}

\title{Coxeter groups of stellar manifolds}

\author{Sergey Nikitin \\
Department of Mathematics\\
and\\
Statistics\\
Arizona State University\\
Tempe, AZ 85287-1804}

\maketitle

{\bf Abstract} It is well known that a compact two dimensional surface is homeomorphic to a polygon with the edges identified in pairs. This paper not only presents a new proof of this statement but also generalizes it to any connected $n$-dimensional stellar manifold with a finite number of vertices. The study of an $n$-dimensional stellar sphere with the boundary faces identified in pairs leads us to Coxeter groups of stellar manifolds. Analysis of these Coxeter groups allows us to single out the class of so called flat stellar manifolds. Flat stellar manifolds are not only very common but also possess a range of simple properties that allow to prove some difficult long standing statements. In particular, there exists a simple proof of the Poincar\`e conjecture for flat manifolds. It is not possible to carry out this simple proof for a general stellar manifolds because of certain singularities that can be characterized in terms of Coxeter groups. The results and ideas presented in this paper not only open a new way that may lead to entirely combinatorial and algebraic proof of the Poincar\`e conjecture but also suggest a new technique for studying stellar manifolds in terms of properties of their Coxeter groups.

\vspace{0.1cm}

\section{Introduction}

\vspace{0.1cm}

It is well known (see e.g. \cite{Massey}) that a compact $2$-dimensional surface is homeomorphic to a polygon with the edges identified in pairs.
 In this paper we show that any connected $n$-dimensional stellar manifold $M$ with a finite number of vertices possesses a similar property: $M$ is stellar equivalent to 
$$
       a\star (S/\simeq),
$$
where $a \notin S$ is a vertex, $S$ is $(n-1)$-dimensional stellar sphere. Throughout this paper it is tacitly assumed that the term "regular equivalence" (denoted as "$\simeq$") actually comprises two equivalence relations: one on the set of vertices and the other on the set of $(n-1)$-dimensional simplexes of $S.$ Those two equivalence relations satisfy the following conditions:
\begin{itemize}
\item[] No $(n-1)$-dimensional simplex of $S$ has two vertices that are equivalent to each other.
\item[] For any $(n-1)$-dimensional simplex $g$ of $S$ there might exist not more than one (different from $g$) $(n-1)$-dimensional simplex $p$ in $S$ such that$p\simeq g$ and any vertex of $g$ is equivalent to some vertex of $p.$
\end{itemize}
$a \star (S/\simeq)$ is not any more a stellar manifold in its usual sense.
We take the ball $a\star S$ and identify equivalent simplexes from $S$ but,
at the same time, we distinguish any two simplexes $a\star p$ and $a\star g$ even when $p\simeq g.$ In the sequel, 
it is convenient to use the notation $a\star (S/\simeq)$ for the resulted manifold even if it is slightly misleading. $a\star (S/\simeq)$ is called a stellar structure of $M.$\\
The complex $S/\simeq$ inherits certain properties from $a\star (S/\simeq).$ Moreover,  $S/\simeq$ can be entirely described  by its Coxeter group.

\section{Preliminaries from stellar theory}
We begin with recalling the basic definitions of stellar theory \cite{Glaser}, \cite{Lickorish}.
A stellar $n$-manifold $M$ can be identified with the sum of its $n$-dimensional simplexes ($n$-simplexes):
$$
M=\sum_{i=1}^n g_i
$$
with coefficients from ${\rm Z}_2.$ We will call $\{g_i\}_{i=1}^n$ generators of $M.$

All vertices in $M$ can be enumerated and any $n$-simplex $s$ from $M$ corresponds to the set
of its vertices
$$
s=(i_1 \; i_2 \; \dots \; i_{n+1}),
$$
where $i_1 \; i_2 \; \dots \; i_{n+1}$ are integers. 

The boundary operator $\partial $ is defined on a simplex as
$$
\partial (i_1 \; i_2 \; \dots \; i_{n+1}) = (i_1 \; i_2 \; \dots \; i_n ) + (i_1 \; i_2 \; \dots \; i_{n-1} \; i_{n+1} ) + \dots + (i_2 \; i_3 \; \dots \; i_{n+1})
$$
and linearly extended to any complex, i.e.
$$
\partial M = \sum_{i=1}^n \partial g_i.
$$
A manifold is called closed if $\partial M =0.$

If two simplexes $(i_1 \; i_2 \; \dots i_m) $ and $(j_1 \; j_2 \; \dots j_n)$ do not have common vertices then one can define their join
$$
(i_1 \; i_2 \; \dots i_m) \star (j_1 \; j_2 \; \dots j_n)
$$
as the union
$$
(i_1 \; i_2 \; \dots i_m) \cup (j_1 \; j_2 \; \dots j_n).
$$

\vspace{0.1cm}

If two complexes $K=\sum_i q_i $ and $ L = \sum_j p_j$ do not have common vertices then 
their join is defined as
$$
K\star L = \sum_{i,j} q_i \star p_j.
$$

If $A$ is a simplex in a complex $K$ then we can introduce its link:
$$
lk(A,K) = \{ B \in K \; ; \; A \star B \in K \}.
$$
The star of $A$ in $K$ is $A \star lk(A,K).$ Thus,
$$
K = A\star lk(A,K) + Q(A,K),
$$
where the complex $Q(A,K)$ is composed of all the generators of $K$ that do not contain $A.$ A complex with generators of the same dimension is called a uniform complex.

\begin{definition} ({\bf Subdivision})
Let $A$ be a simplex of a complex $K.$ Then any integer $a$ which is not a vertex of $K$ defines starring of
$$
K=A\star lk(A,K) + Q(A,K)
$$
at $a$ as
$$
 \hat{K}=a\star \partial A \star lk(A,K) + Q(A,K).
$$
This is denoted as
$$
\hat{K} =  (A \; a)K.
$$
\end{definition}

The next operation is the inverse of subdivision. It is called a stellar weld and defined as follows.

\begin{definition} ({\bf Weld})
  Consider a complex
  $$
        \hat{K}=a\star  lk(a,\hat{K}) + Q(a,\hat{K}),
  $$
with $lk(a,\hat{K}) = \partial A \star B$ where
$B$ is a subcomplex in $\hat{K},$  $A$ is a simplex and $A\notin \hat{K} .$
Then the (stellar) weld $(A\; a)^{-1} \hat{K}$ is defined as
  $$
       (A\; a)^{-1} \hat{K} =  A \star B  + Q(a,\hat{K}).
  $$
\end{definition}

A stellar move is one of the following operations: subdivision, weld, enumeration change on the set of vertices.
Two complexes $M$ and $L$  are called stellar equivalent if one is obtained from the other by a finite sequence of stellar moves.
It is denoted as $M \sim L.$ We also say that $M$ admits triangulation $L.$\\

If a complex $L$ is stellar equivalent to $(1 \; 2 \; \dots \; n+1)$ then $L$ is called a stellar 
$n$-ball. On the other hand, if $K \sim \partial (1 \; 2 \; \dots \; n+2)$ then $K$ is a stellar
$n$-sphere.

\begin{definition} ({\bf Stellar manifold})
\label{stellar_def}
Let $M$ be a complex. If, for every vertex $i$ of $M,$ the link $lk(i,M)$ is either a stellar $(n-1)$-ball or a stellar $(n-1)$-sphere, then $M$ is  a
stellar $n$-dimensional manifold ($n$-manifold).
\end{definition}

If $i$ is a vertex of $M$ then
$$
M= i\star lk(i,M) + Q(i,M).
$$
If $\partial M=0,$ then $Q(i,M)$ is a stellar manifold.\\
Indeed, consider an arbitrary vertex $j$ of $Q(i,M).$ Then
$$
lk(i,M)=j\star lk(j,lk(i,M)) + Q(j,lk(i,M))
$$
and
$$
Q(i,M)=j\star lk(j,Q(i,M)) + Q(j,Q(i,M)).
$$
Since $M$ is a stellar manifold and $\partial M = 0$
$$
i\star lk(j,lk(i,M)) + lk(j,Q(i,M)) 
$$
is a stellar sphere. Hence, it follows from \cite{Newman} that $lk(j,Q(i,M))$ is either a stellar ball or a stellar sphere.\\

The following theorem plays one of the central roles in the research on stellar manifolds.  

\begin{theorem} (Alexander \cite{Alexander})
\label{Alexander}
Let $M$ be a stellar $n$-manifold, let $J$ be a stellar $n$-ball. Suppose that
$M\cap J = \partial M \cap \partial J$ and that this intersection is a stellar $(n-1)$-ball. Then $M\cup J$ is stellar equivalent to $M.$
\end{theorem}

One can say that $M$ is obtained from $M\cup J$ by collapsing a stellar $n$-ball. If a  stellar $n$-manifold $L$ is obtained from  a stellar $n$-manifold $M$ by collapsing a finite number of simplexes (not necessary having dimension $n)$ then $M$ collapses to $L.$  If $L$ is a point then $M$ is called collapsible. \\
It was shown by Whitehead \cite{Whitehead} that a collapsible stellar manifold is a stellar ball. This fact is of special interest to us because it implies the Poincar\'e conjecture for flat stellar manifolds. On the other hand, any stellar manifold can be constructed out of finite number of flat pieces that can be glued together along singularities classified with the help of finite Coxeter groups. These array of facts opens a way of dealing with Poincar\'e conjecture on the basis of algebraic and combinatorial technique.\\ 

For reader convenience we present here a definition of prism over a complex \cite{Alexandrov}.
Let $K$ be a uniform $n$-dimensional complex with vertices $i_1, \; i_2,  \; \dots \;  i_n .$ Let us take new vertices $j_1,  \; j_2 ,\; \dots \; j_n$ which correspond one-to-one to 
the vertices $i_1 , \; i_2, \; \dots \; i_n .$
The complex $P(K)$ is called a prism over $K$ if $P(K)$ consists, by definition, of all nonempty subsets of sets of the form
$$
 i_1 ,\; i_2 , \; \dots  ,\; i_k, \; j_k ,\; \dots ,\;j_r,
$$
where 
$$
 (i_1 ,\; i_2 , \; \dots  ,\; i_r) \in K.
$$
Moreover, $P(K)$ contains all the simplexes $(i_1 ,\; i_2 , \; \dots  ,\; i_r)$ from $K$ and all the simplexes $(j_1 ,\; j_2 , \; \dots  ,\; j_r)$ corresponding to them.\\
The prism $P(K)$ can be also constructed in the following way. Consider 
$$
a \star K,
$$
where $a\notin K$ is a vertex. Take its subdivision defined as
$$
L=\prod_{b \in K} ((a\; b)\; c_b) K,
$$
where $\prod_{b \in K} ((a\; b)\; c_b)$ is a superposition of subdivisions $((a\; b)\; c_b),$ where $c_b \notin K$ and $c_b \not= c_d $ for $b \not= d.$ Then
$$
L = a \star lk(a,L) + Q(a,L),
$$
where $Q(a,L)=P(K)$ is a prism over $K$ and $lk(a,L)$ is obtained from $K$ by the enumeration change $b \rightarrow c_b$ for $b\in K.$\\

In the sequel it is convenient to consider two equivalence relations: one on the set of vertices  and the other on the set of generators of a stellar manifold. Among all possible such equivalence relations we are mostly interested in those that meet certain regularity properties underlined by the following definition.
\begin{definition} ({\bf Regular equivalence})
Given a stellar manifold $M,$ a pair of equivalence relations,  one on the set of vertices and the other on the set of generators from $M,$ is called regular equivalence if it meets the following conditions:
\begin{itemize}
\item[(i)] No generator $g \in M$ has two vertices that are equivalent to each other.
\item[(ii)]For any generator $g\in M$ there might exist not more than one generator $p\in M $ such that $p \not= g$ and $p$ is equivalent to $g,\;\;g \simeq p.$ Moreover, any vertex of $g$ is either equal or equivalent to some vertex of $p.$ 
\end{itemize}
\end{definition}

Throughout the paper $a\star (S/\simeq)$ denotes the structure obtained from $a\star S$ when we identify any two equivalent generators from $S$ but distinguish the corresponding generators from $a\star S.$

\section{Stellar Manifold Structure }
We begin with proving that any connected stellar $n$-manifold with a finite number of generators can be obtained from a stellar $n$-ball $B$ by identifying in pairs the generators of $\partial B.$ For $3$-manifolds this fact was stated in \cite{SeifThrelf}.

\begin{theorem}({\bf Stellar Structure})
\label{triangulation}
A connected stellar $n$-manifold $M$ with a finite number of generators admits a triangulation (stellar structure)
$$
N = a\star (S/\simeq) ,
$$
where $a\notin S$ is a vertex, $S$ is a stellar $(n-1)$-sphere and "$\simeq $" is a regular equivalence relation. Moreover, if $M$ is closed then for any generator $g\in S$ there exists exactly one generator $p\in S$ such that $p\not= g$ and $g \simeq p.$ 
\end{theorem}
{\bf Proof.}
Let us choose an arbitrary generator $g \in M$ and an integer $a$ that is not a vertex of $M.$ Then
$$
M \sim (g \; a ) M  \mbox{  and  }
(g \; a ) M =a \star \partial g + M\setminus g,
$$
where $M\setminus g$ is defined by all the generators of $M$ excluding $g.$
We construct $N$ in a finite number of steps.
Let $N_0=(g \; a ) M.$ Suppose we constructed already $N_k$ and there exists a generator $p\in  Q(a,N_k)$ that has at least one common $(n-1)$-simplex with $lk(a, N_k).$
Without loss of generality, we can assume that 
$$
 p = (1 \; 2\; \dots \; n+1).
$$
and $(1 \; 2\; \dots \; n)$ belongs to $lk(a, N_k).$
If the vertex $(n+1)$ does not belong to  $lk(a, N_k)$ then
$$
N_{k+1} = ((a\; n+1 )\; b)^{-1}((1 \; 2 \; \dots \; n) \; b) N_k,
$$
where $b\notin N_k.$
If the vertex $(n+1)$ belongs to $lk(a, N_k)$ then after introducing a new vertex $d\notin N_k$ we take
$$
L = ((a\; d )\; b)^{-1}((1 \; 2 \; \dots \; n) \; b)(N_k\setminus p + (1\;2\;\dots \; n \;d)),
$$
where $b\notin (N_k\setminus p + (1\;2\;\dots \; n \;d )),$ and
$$
N_{k+1} =  a \star (lk(a,L)/\simeq ) + Q(a,L)
$$
endowed with the equivalence $d \simeq (n+1).$\\
By construction
$$
N_{k+1} = a \star (lk(a,N_k) +\partial p ) + Q(a,N_k) \setminus p
$$
if $(n+1)\notin lk(a, N_k).$ Otherwise,
$$
N_{k+1} = a \star ((lk(a,N_k) +\partial g)/\simeq ) + Q(a,N_k) \setminus p,
$$
where $g=(1 \; 2\; \dots \; n \; d),\;\;d\simeq (n+1).$

Since $M$ is connected and has a finite number of generators there exists a natural number $m$ such that
$$
N_m =  a \star (S/\simeq)
$$
where $S$ is a stellar $(n-1)$-sphere and "$\simeq $" is a regular equivalence relation.\\
If $M$ is closed, then $\partial N_m = 0,$ and therefore, for any generator $g \in S$ there exists exactly one generator $p\in S \setminus g$ such that $g \simeq p.$\\
Q.E.D.\\

When working with stellar manifolds we use topological terminology in the following sense. We say that a stellar $n$-manifold $M$ has a certain topological property if its standard realization in ${\rm R}^{2n+1}$ (see e.g \cite{Alexandrov}) has this property.\\
Let $\pi (M)$ denote the fundamental group of a manifold $M.$
Then Theorem \ref{triangulation} together with Seifert -- Van Kampen theorem (see e.g. \cite{Massey}) lead us to the next result.

\begin{theorem}
\label{fundamental}
If $M \sim a \star (S/\simeq )$ is a connected stellar $n$-manifold and $n>2$ then
$$
\pi (M) = \pi (S/\simeq ),
$$
where $S$ is a stellar $(n-1)$-sphere and "$\simeq $" is a regular equivalence relation.
\end{theorem}
{\bf Proof.}
By Theorem \ref{triangulation}  $M$ admits a triangulation
$$
N = a\star (S/\simeq ),
$$
where $a\notin S$ is a vertex, $S$ is a stellar $(n-1)$-sphere and "$\simeq $" is a regular equivalence relation.\\
Using the triangulation $N$ we can define open subsets $U$ and $V$ in $M$ as follows
\begin{eqnarray*}
U &=& a\star (S/\simeq )\setminus a \\
V &=& a\star (S/\simeq ) \setminus (S/\simeq ).
\end{eqnarray*}
Clearly,
$$
M = U \cup V
$$
and $V$ is homeomorphic to an interior of the $n$-ball. Thus, $V$ is simply connected and by Seifert -- Van Kampen theorem there exists an epimorphism
$$
\psi \;:\; \pi(U)\; \longrightarrow \; \pi(M)
$$
induced by inclusion $U \subset M .$ The kernel of $\psi $ is the smallest normal subgroup containing the image of the homomorphism
$$
\varphi\;:\;\pi(U\cap V) \; \longrightarrow \; \pi(U)
$$
induced by inclusion $U\cap V \subset U.$\\

The open set $U\cap V $ is homeomorphic to 
$$
0 < x_1^2+ x_2^2 +\dots +x_n^2 <1,\;\;\mbox{ where } \;\; (x_1,\; x_2 ,\; \dots ,\;x_n) \in {\rm R}^n.
$$
The fundamental subgroup of this manifold is trivial for $n>2.$ Hence, the kernel of the epimorphism $\psi $ is trivial and
$$
\pi(M) = \pi(U).
$$
On the other hand, $S/\simeq $ is a deformation retract of $U,$ and therefore,
$$
 \pi(U) = \pi (S/\simeq).
$$
Q.E.D.\\

Let $\chi (M)$ denote Euler characteristic of $M,$ i.e. 
$$
\chi(M) = \sum_{i=0}^n (-1)^i q_i,
$$
where $q_i$ is the number of $i$-simplexes in $M.$

\begin{theorem}
\label{Euler}
If $M \sim a \star (S/\simeq )$ is a closed connected stellar $n$-manifold then 
$$
      \chi(S/\simeq ) = \chi(M) + (-1)^{n+1},
$$
where $S$ is a stellar $(n-1)$-sphere and "$\simeq $" is a regular equivalence relation.
\end{theorem}
{\bf Proof.}
Since $S$ is a stellar $(n-1)$-sphere we have
$$
\chi (S) = \sum_{i=0}^{n-1} (-1)^i s_i = (-1)^{n-1}+1,
$$
where $s_i$ is the number of $i$-simplexes in $S.$
On the other hand,
$$
\chi( a \star (S/\simeq ))=\sum_{i=0}^n (-1)^i q_i,
$$
where $q_i$ is the number of $i$-simplexes in $a \star (S/\simeq ) .$ 
Let $h_i$ denote the number of $i$-simplexes in $S/\simeq .$
Taking into account that "$\simeq$" is a regular equivalence relation we obtain
$$
s_{n-1} = 2 h_{n-1},\;\;q_n = 2 h_{n-1} \;\; \mbox{ and } \;\; q_i = h_i + s_{i-1}\;\;\mbox{ for }\;\; 1\le i \le n-1
$$
and
$$
q_0 = h_0 + 1.
$$
Thus,
$$
\chi(M) = \sum_{i=0}^n (-1)^i q_i= 2 (-1)^n h_{n-1}  + h_0 + 1+ \sum_{i=1}^{n-1} (-1)^i(h_i + s_{i-1})
$$
and
$$
 \sum_{i=1}^{n-1} (-1)^i s_{i-1} = (-1)^{n-1} s_{n-1}+ (-1)^n-1 = (-1)^{n-1} 2 h_{n-1} + (-1)^n-1.
$$
Hence,
$$
\chi(M) =  \chi(S/\simeq ) + (-1)^n
$$
and the assertion is proved.\\
Q.E.D.\\

\section{Coxeter group of a stellar structure}
This section deals only with finite Coxeter groups \cite{Bourbaki}, \cite{Humphreys}. A Coxeter group is a group with a finite set of generators $\Pi$  such that
$$
\forall r \in \Pi \;\;\; r^2 =1 \mbox{  and  } \forall r,\;q \in  \Pi\;\; \exists\; m_{rq} \mbox{ such that } (r\cdot q )^{m_{rq}}=1,
$$
where $m_{rr} = 1$ and $ m_{rq} \ge 2 $ for $r \not= q.$
A stellar structure $a\star (S/\simeq)$ possesses a Coxeter group defined as follows.   

\begin{definition}
\label{CoxeterGroup}
 Consider a closed stellar manifold $a\star (S/\simeq).$ Let ${\cal G}=\{g_1,\;\dots ,g_m  \}$ be the set of generators of the stellar $n$-sphere $S.$ The regular equivalence $\simeq$ defines the permutation $p_0$ on $\cal G,$ 
$$
  p_0(g_i) = g_j \mbox{ for } g_i \simeq g_j \mbox{ and } g_i \not= g_j.
$$
We call two $(n-1)$-simplexes $\alpha $ and $\beta $ equivalent if they belong to equivalent $n$-simplexes in $S$ and each vertex of  $\alpha $ is either equal or equivalent to some vertex in $\beta .$ Let $\cal E$ denote all equivalence classes of $(n-1)$-simplexes in $S$.  Then each $\alpha \in {\cal E}$ corresponds to a permutation $p_{\alpha }$ defined as follows. If $g_i \in {\cal G}$ does not have any $(n-1)$-simplexes from $\alpha$ then
$$
    p_\alpha (g_i) = g_i.
$$
Otherwise, $g_i$ may have only one $(n-1)$-simplex from $\alpha.$ There exists only one $g_j \in  {\cal G}$ that contains the same $(n-1)$-simplex. We define 
$$
   p_\alpha (g_i) = g_j.
$$
The Coxeter group of $S/\simeq$ is defined by the generators $\{p_0,\;\{p_\alpha \}_{\alpha \in {\cal E}}\}.$
\end{definition}

In the sequel $Cox(S/\simeq)$ denotes the Coxeter group of the stellar structure $a \star (S/\simeq).$ 

\begin{theorem}
\label{flatTheorem}
 Let $M=a\star (S/\simeq )$ be an $(n+1)$-dimensional closed stellar manifold. Then $Cox(S/\simeq )$ has the following property.
\begin{equation}
\label{inv}
        (p_0 \cdot p_\alpha )^2 = 1
\end{equation} 
if, and only if, the equivalence class $\alpha $ contains not more than two different $(n-1)$-simplexes. 
\end{theorem}

{\bf Proof.}
 If  $\alpha $ contains a single $(n-1)$-simplex $i=(i_1,\dots , i_{n-1})$ then there exist two vertices $\nu$ and $\mu$ such that
\begin{eqnarray*}
         p_0(\nu \star i) &=& \mu \star i\\
         p_\alpha (\nu \star i) &=& \mu \star i
\end{eqnarray*}
 If $g$ is an $n$-simplex different from $\nu \star i$ and $ \mu \star i$ then $p_\alpha(g)=g$ and
$$
          (p_0\cdot p_\alpha)^2(g)= p^2_0(g)=g.
$$ 
On the other hand,
$$
 (p_0 \cdot p_\alpha)(\nu \star i)= \nu \star i,
$$
and
$$
 (p_0 \cdot p_\alpha)(\mu \star i)= \mu \star i.
$$
Hence, we established (\ref{inv}) when $\alpha $ contains a single $(n-1)$-simplex. \\

If $\alpha $ contains two different $(n-1)$-simplexes $i$ and $j$ then there exist vertices $v_i$, $w_i$, $v_j$ and $w_j$ such that for the following generators of $S$ 
\begin{equation}
\label{eq_class}
v_i\star i,\; w_i \star i,\; v_j \star j,\; w_j\star j 
\end{equation}
we have
\begin{eqnarray*}
 p_\alpha (v_i\star i) &=& w_i \star i \\
 p_\alpha (v_j\star j) &=& w_j \star j \\
 p_0 (v_j\star j) &=& v_i\star i.
\end{eqnarray*}
Since $\alpha $ has only two elements we obtain
$$
p_0 (w_j\star j) = w_i\star i.
$$
If a generator $g$ is not one of the simplexes from (\ref{eq_class}) then 
$p_\alpha (g)=g,$ and $(p_0\cdot p_\alpha )^2 (g)=g.$ 
On the other hand,
$$
(p_0\cdot p_\alpha )^2 (v_i\star i)= p_0\cdot p_\alpha \cdot p_0 (w_i \star i)=p_0\cdot p_\alpha (w_j\star j) = p_0 (v_j\star j)= v_i\star i.
$$
Similar calculation show that
\begin{eqnarray*}
(p_0\cdot p_\alpha )^2 (w_j\star j) & =& w_j\star j \\
(p_0\cdot p_\alpha )^2 (v_j\star j) & =& v_j\star j \\
(p_0\cdot p_\alpha )^2 (w_i\star i) & =& w_i\star i .
\end{eqnarray*}
Hence, (\ref{inv}) is established when $\alpha $ contains not more than two different simplexes.\\

Now we need to show that (\ref{inv}) implies that the equivalence class $\alpha$ contains not more than two different $(n-1)$-simplex. We conduct the proof by reductio ad absurdum. Suppose  $\alpha$ contains at least three different  $(n-1)$-simplexes: $i,$ $j$ and $k.$ Then there exist three vertices $\mu, \;\;\nu $ and $\gamma$ such that $\mu \star i, \;\; \nu \star j$ and $\gamma \star k$ are generators from $S$ and
\begin{eqnarray*}
 p_0(\mu \star i)&=& \nu \star j\\
 p_0\cdot p_\alpha (\nu \star j)&=& \gamma \star k.
\end{eqnarray*}
Hence,
$$
 (p_\alpha \cdot p_0)^2(\mu \star i) = p_\alpha \cdot p_0 \cdot p_\alpha( \nu \star j)= p_\alpha (\gamma \star k).
$$
However,
$$
p_\alpha (\gamma \star k) \not= \mu \star i.
$$
Thus, 
$$
 (p_\alpha \cdot p_0)^2 \not= 1.
$$

Q.E.D.\\

If $(p_0 \cdot p_\alpha )^2=1$ then $(n-1)$-simplex $\alpha$ may belong to a single $n$-simplex or can be shared by two $n$-simplexes of $S/\simeq.$  We underline this fact by calling $S/\simeq$ flat at $\alpha .$

\begin{theorem}
\label{order2}
Let $M=a\star (S/\simeq)$ be a closed $(n+1)$-dimensional stellar manifold.
If $\alpha$ and $\beta$ are two different $(n-1)$-simplexes from $S/\simeq$ then the following is true.  
\begin{equation}
\label{inv1}
  (p_\alpha \cdot p_\beta )^2 =1
\end{equation}
if, and only if, one of the following statements take place.\\
\begin{itemize}
\item[(i)] There is no $n$-simplex in $S$  that contains $(n-1)$-simplexes from both equivalence classes $\alpha$ and $\beta.$
\item[(ii)]If there is an  $n$-simplex $\gamma \in S$ that contains $(n-1)$-simplexes from both equivalence classes $\alpha$ and $\beta$  then this simplex is of the form 
$$
\gamma = i\star j \star \ell,
$$
where $i,\; j$ are vertices of $S;$  $\ell $ is an $(n-2)$-simplex in $S$ such that
\begin{equation}
\label{starell}
lk( \ell, S) = (i\; j) + (j\; k) + (k\; m ) + (m \; i),
\end{equation}
where $k,\;m$ are vertices in $S$ and
\begin{equation}
\label{stareq}
 (i \star \ell ) \simeq (k \star \ell ), \;\;\; (m \star \ell ) \simeq ( j \star \ell ).
\end{equation}
\end{itemize}
\end{theorem}
{\bf Proof.}
If an $n$-simplex $g\in S$ does not contain any $(n-1)$-simplex from $\alpha $ then
$$
 p_\alpha (g) = g.
$$
Likewise, 
$$
p_\beta (q) = q
$$
if an $n$-simplex $q$ does not contain any $(n-1)$-simplex from $\beta .$
Thus, statement (i) implies (\ref{inv1}).\\
On the other hand, if an $n$-simplex $\gamma $  contains $(n-1)$-simplexes from both equivalence classes $\alpha$ and
$\beta $ then 
$$
\gamma = i\star j \star \ell,
$$
where $j\star \ell$ is a representative of $\alpha$ and  $i\star \ell$ belongs to $\beta.$
$$
p_\alpha (i\star j \star \ell) = p_\alpha (i)\star j \star \ell,
$$
where $p_\alpha (i)$ the corresponding vertex of $p_\alpha (i\star j \star \ell).$ 
It follows from (\ref{inv1}) that 
$$
(p_\beta \cdot p_\alpha )^2 (p_\alpha (i)\star j \star \ell) = p_\alpha (i)\star j \star \ell .
$$
That means
$$
 p_\beta \cdot p_\alpha \cdot p_\beta ( i\star j \star \ell)=  p_\beta \cdot p_\alpha (i\star p_\beta(j) \star \ell).
$$
If $p_\beta(j) \star \ell$ is not in $\alpha $ then 
$$
  p_\beta \cdot p_\alpha (i\star p_\beta(j) \star \ell) = p_\beta (i\star p_\beta(j) \star \ell) = i\star j \star \ell
$$
and (\ref{inv1}) is not valid. 
Hence, $p_\beta(j) \star \ell$ is equivalent to $j\star \ell$ and
$$
  p_\beta \cdot p_\alpha (i\star p_\beta(j) \star \ell) = p_\beta (\bar p_\alpha (i) \star p_\beta(j) \star \ell).
$$
It follows from (\ref{inv1}) that $\bar p_\alpha (i) \star \ell$ is in $\beta $ and
$$
p_\beta (\bar p_\alpha (i) \star p_\beta(j) \star \ell) = \bar p_\alpha (i) \star p^2_\beta(j) \star \ell = p_\alpha (i)\star j \star \ell.
$$
Therefore, $\bar p_\alpha (i) = p_\alpha (i),$  $p^2_\beta(j) = j.$ After introducing notations $p_\alpha (i) = k, \; p_\beta(j)=m$ we obtain (\ref{starell}), (\ref{stareq}). The sufficiency of conditions (i) and (ii) is evident.
The proof is completed.
Q.E.D.\\

Now we analyze geometrical properties of $S/\simeq$ imposed by algebraic conditions
\begin{eqnarray*}
 (p_0 \cdot p_\alpha )^m &=& 1 \mbox{ and }  (p_0 \cdot p_\alpha )^{s} \not= 1 \mbox{ for } 0<s < m\\
 (p_\alpha \cdot p_\beta)^k &=& 1  \mbox{ and }  (p_\alpha \cdot p_\beta)^{s} \not= 1 \mbox{ for } 0<s < k,
\end{eqnarray*}
where natural numbers $k,m$ are larger than $2.$

\begin{theorem}
\label{orderM}
Let $M=a\star (S/\simeq )$ be a closed $(n+1)$-dimensional stellar manifold. Then
\begin{equation}
\label{power}
 (p_0 \cdot p_\alpha )^m = 1  \mbox{ and }  (p_0 \cdot p_\alpha )^{k} \not= 1 \mbox{ for } 0<k < m
\end{equation}
for some integer $m>2$ if, and only if, the equivalence class $\alpha$ contains $m$ different $(n-1)$-simplexes from $S.$
\end{theorem}

{\bf Proof.}
Let (\ref{power}) be valid. If $j$ is an $(n-1)$-simplex from the equivalence class $\alpha $ then there exists a vertex $\nu$ from $S$ such that $\nu \star j$ is an $n$-simplex from $S$ and
$$
 p_\alpha (\nu \star j) = p_\alpha (\nu ) \star j,
$$
where $p_\alpha (\nu ) $ is the corresponding vertex from $S.$
Since $m>2$ $p_0 (p_\alpha (\nu ) \star j) $ is different from $ \nu \star j.$ Introduce notations
$$
 p_0 (p_\alpha (\nu ) \star j) = p_0\cdot p_\alpha (\nu) \star p_0 (j),
$$
where $p_0\cdot p_\alpha (\nu)$ is a vertex and $p_0 (j)$ is a $(n-1)$-simplex from $S.$
The equivalence class $\alpha $ contains $p_0 (j)$, and therefore, 
$$
p_\alpha (p_0\cdot p_\alpha (\nu) \star p_0 (j))= p_\alpha \cdot p_0\cdot p_\alpha (\nu )  \star p_0 (j)
$$
is well defined.
The $n$-simplex 
$$
p_0(p_\alpha \cdot p_0\cdot p_\alpha (\nu )  \star p_0 (j))= (p_0\cdot p_\alpha)^2 (\nu) \star p^2_0(j)
$$
is different from $\nu \star j$ because $m>2.$
If $m=3$ then
$$
 (p_0\cdot p_\alpha)^3(\nu) \star p^3(j) = \nu \star j,
$$
where $(p_0\cdot p_\alpha)^3 (\nu)$ is a vertex  and $ p^3(j)$ is the corresponding $(n-1)$-simplex of $(p_0\cdot p_\alpha)^3(\nu \star j).$
If $m > 3$ then we continue in this fashion and conclude that
$$
(p_0\cdot p_\alpha)^k \star p^k(j) = \nu \star j
$$
if, and only if, $k=m.$ Hence, $\alpha $ contains $m$ different $(n-1)$-simplexes
$$
j,\; p_0(j),\; \dots , p^{m-1}_0 (j).
$$
For any $n$-simplex g from $S$ that does not contain an $(n-1)$-simplex from $\alpha$ we have 
$$
 (p_0 \cdot p_\alpha )^2 (g) = g.
$$
If $\alpha $ contains $m$ different $(n-1)$-simplexes
$$
j,\; p_0(j),\; \dots , p^{m-1}_0 (j).
$$
then there exits a vertex $\nu$ such that
$$
 (p_0 \cdot p_\alpha )^m (\nu \star j) = \nu \star j
$$
and
$$
 (p_0 \cdot p_\alpha )^k (\nu \star j) \not= \nu \star j \mbox{ for }  0<k<m.
$$
Q.E.D.\\

\begin{figure}[tb]
  \begin{center}
     \psfig{file=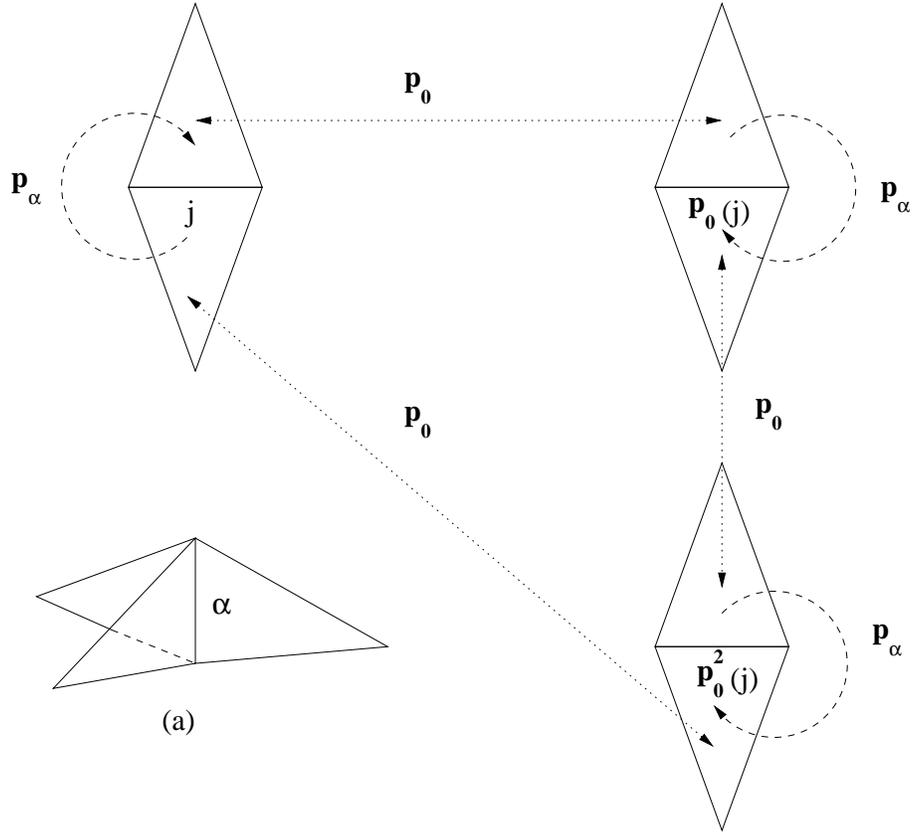,scale=0.7}         
      \caption{$(p_0 \cdot p_\alpha )^3 =1;$ (a) -- a view of $S/\simeq $ near $\alpha.$}
    \label{fig:order3} 
  \end{center}
\end{figure}

The statement of Theorem \ref{orderM} is illustrated in Fig.\ref{fig:order3} for $m=3.$ The geometrical appearance of $S/\simeq $ is sketched in Fig.\ref{fig:order3} (a).

\begin{definition}
Let $M=a\star (S/\simeq)$ be a closed $(n+1)$-dimensional stellar manifold. Given two $(n-1)$-simplexes $\alpha$ and $\beta$ in $S/\simeq$ the order of an $n$-simplex $g\in S $ (with respect to $\alpha $ and $\beta $) is called the smallest integer $m>0$ such that
$$
  (p_\alpha \cdot p_\beta )^m (g) =g.
$$
The order of $g$ (with respect to $\alpha $ and $\beta $) is denoted by 
$$
order_{\alpha  \beta } (g).
$$
\end{definition}
The order is symmetric,
$$
   order_{\alpha  \beta } (g) = order_{\beta  \alpha } (g) \;\;\;\forall\; g.
$$
If $order_{\alpha  \beta } (g) > 2 $ then $g$ is of the form
\begin{equation}
\label{powerS}
g = i\star j \star \ell,
\end{equation}
where $i,\; j$ are vertices of $S$ and $\ell $ is an $(n-2)$-simplex in $S$ such that $j\star \ell$ is in $\beta $ and $i \star \ell $ is in $\alpha.$ Moreover, the inequality 
$$
     order_{\alpha , \beta } (g) > 2 
$$
dictates for the star $\ell \star lk(\ell,S)$ to have a certain structure described in the next theorem.
\begin{theorem}
\label{morethan2}
 Let $a\star (S/\simeq)$ be a stellar structure of a closed $(n+1)$-dimensional stellar manifold. If  $\alpha$ and $\beta$ are two $(n-1)$-simplexes in $S/\simeq$ and $g$ is an $n$-simplex in $S$ then
\begin{equation}
\label{power1}
   order_{\alpha  \beta } (g)  >2
\end{equation}
if, and only, if  $g$ takes form (\ref{powerS}) and the star  $\ell \star lk(\ell,S)$ satisfies one of the following properties.

\begin{itemize}
\item[(i)] If $order_{\alpha  \beta } (g) = m$ then there exist exactly $2\cdot m$ vertices $\{i_1,\;i_2, \dots , i_{2m} \}$ such that
$$
 lk(\ell,S) = (i_1 \;i_2) + (i_2 \; i_3)+ \dots + (i_{2m - 1}\; i_{2m}) + (i_{2m} \; i_1),
$$
and $i_{2k-1} \star \ell \in \alpha , \; i_{2k} \star \ell \in \beta $ for $k=1,\;2, \dots , m.$
\item[(ii)] There exists a vertex $v \in lk(\ell,S)$ such that $v \star \ell$ is neither in $\alpha $ nor in $\beta .$ Moreover, if  $r_{\alpha \beta },\;\;r_{\beta \alpha }$ are the smallest natural numbers for which 
\begin{eqnarray}
\label{ralphabeta}
 (p_\alpha \cdot p_\beta )^{r_{\alpha \beta }} (g) &=& p_\beta (g) \nonumber \\
 && \\
 (p_\beta \cdot p_\alpha )^{r_{\beta \alpha }} (g) &=& p_\alpha (g) \nonumber
\end{eqnarray}
then
$$
     order_{\alpha  \beta } (g) = r_{\alpha \beta} + r_{\beta \alpha } + 1.
$$
\end{itemize}
\end{theorem}
{\bf Proof.}
(\ref{power1}) implies the existence of an $n$-simplex 
$$
g= i\star j \star \ell,
$$
where $j \star \ell$ is in $\beta $ and $ i\star \ell$ is in $\alpha .$ Otherwise, by Theorem \ref{order2},
$$
 (p_\alpha \cdot p_\beta )^2=1.
$$
For $lk(\ell, S)$ there are only two possibilities:
\begin{itemize}
 \item[(p1)] For any vertex  $v \in lk(\ell, S)$ the $(n-1)$ - simplex $v \star \ell$ is either in $\beta $ or in $\alpha .$ 
 \item[(p2)] There is a vertex $v \in lk(\ell, S)$ such that $v \star \ell$ is neither in $ \beta $ nor in $\alpha .$
\end{itemize}

Consider (p1). Let $N$ denote the total number of vertices in $lk(\ell, S).$ Since $\alpha $ and $\beta $ are classes of equivalence in $S/\simeq $ and "$ \simeq $" is a regular equivalence, then  none of the $n$-simplexes in $S$ can contain two different $(n-1)$-simplexes from $\beta.$ Hence, one can enumerate vertices $ lk(\ell, S)$ so that
two subsequent vertices belong to the same $n$-simplex, all vertices with odd indices are such that  $i_{2k-1}\star \ell \in \alpha $  while for all vertices with even indices $i_{2k} \star \ell \in \beta .$ Since
\begin{eqnarray*}
(p_\alpha \cdot p_\beta ) (i_s) &=& i_{s+2} \mbox{ for } s=1,\; 2, \dots,\;N-2 \\
(p_\alpha \cdot p_\beta ) (i_{N - 1}) &=& i_1 \\
(p_\alpha \cdot p_\beta ) (i_{N }) &=& i_2 
\end{eqnarray*}
and
$$
(p_\alpha \cdot p_\beta )^m (i) = i\;\;\forall\; i \in lk(\ell, S)
$$
we conclude that $N=2m$ and the necessity of statement (i) is established.\\
Statement (i) also is sufficient for $order_{\alpha \beta }(g) = m$ because
$$
 (p_\alpha \cdot p_\beta ) ((i_{2k - 1}\; i_{2k})\star \ell) = (i_{2k + 1}\; i_{2(k+1)})\star \ell \;\; \mbox{ for } k=1,\;2, \dots , \;m-1
$$
and
$$
(p_\alpha \cdot p_\beta ) ((i_{2m - 1}\; i_{2m})\star \ell) = (i_1\; i_2) \star \ell. 
$$

Consider (p2). Let $v$ be a vertex from statement (ii). Then for an $n$-simplex $s$ that contains $v \star \ell$ we have either $p_\alpha (s) =s $ or $p_\beta (s) = s.$ That implies the existence of natural numbers $r_{\alpha \beta }$ and $r_{\beta \alpha }$ such that (\ref{ralphabeta}) is valid. After applying to the first equation $p_\beta $ and to the second $p_\alpha $ we obtain
\begin{eqnarray*}
 p_\beta \cdot (p_\alpha \cdot p_\beta )^{r_{\alpha \beta }} (g) &=& g \nonumber \\
 p_\alpha \cdot (p_\beta \cdot p_\alpha )^{r_{\beta \alpha }} (g) &=& g \nonumber
\end{eqnarray*}
Hence,
$$
  p_\beta \cdot (p_\alpha \cdot p_\beta )^{r_{\alpha \beta }} (g) = p_\alpha \cdot (p_\beta \cdot p_\alpha )^{r_{\beta \alpha }}
$$
which is equivalent to
$$
    (p_\alpha \cdot p_\beta )^{r_{\alpha \beta } + r_{\beta \alpha } + 1 } (g) = g.
$$
That means $order_{\alpha  \beta } (g)$ is a divisor of $r_{\alpha \beta } + r_{\beta \alpha } + 1,$ 
$$
     order_{\alpha  \beta } (g) | r_{\alpha \beta } + r_{\beta \alpha } + 1 . 
$$
On the other hand, the existence of vertex $v$ from (ii) implies that one can find the smallest natural numbers $k_+$ and $s_+$ such that
$$
   (p_\alpha \cdot p_\beta )^{k_+}(g)=  p_\beta \cdot (p_\alpha \cdot p_\beta )^{s_+} (g).
$$
Thus,
$$
     (p_\alpha \cdot p_\beta )^{k_+}(g)= (p_\beta \cdot _\alpha)^s{_+} p_\beta (g)
$$
which implies
\begin{equation}
\label{lst}
 (p_\alpha \cdot p_\beta )^{k_+ +  s_+ }(g) = p_\beta (g).
\end{equation}
Hence,
$$
 r_{\alpha \beta } = k_+ +  s_+.
$$
Now, let $k_-$ and $s_-$ are the smallest natural numbers for which
$$
(p_\alpha \cdot p_\beta )^{k_-} p_\alpha(g) = p_\beta (p_\alpha \cdot p_\beta )^{s_-} p_\alpha(g) = (p_\beta \cdot p_\alpha)^{s_- + 1} (g)
$$
Then
\begin{equation}
\label{lst1}
(p_\beta \cdot p_\alpha )^{k_- + s_- + 1} (g) =  p_\alpha(g)
\end{equation}
and
$$
 r_{\beta \alpha } = k_- + s_- + 1.
$$
Applying $p_\alpha \cdot p_\beta $ to both parts in (\ref{lst}) we obtain
\begin{equation}
\label{lst2}
 (p_\alpha \cdot p_\beta )^{k_+ +  s_+ + 1} (g) = p_\alpha(g).
\end{equation}
Combining together (\ref{lst1}) and (\ref{lst2}) we have 
$$
 (p_\alpha \cdot p_\beta )^{k_+ +  s_+ + 1} (g) = (p_\beta \cdot p_\alpha )^{k_- + s_- + 1} (g)
$$
Thus,
$$
 (p_\alpha \cdot p_\beta )^{k_+ +  s_+ + k_- + s_- + 2} (g) = g.
$$
Hence, $k_+ +  s_+ + k_- + s_- + 2$ is a divisor of $order_{\alpha  \beta } (g)$ and so is
$
 r_{\alpha \beta } + r_{\beta \alpha } + 1.
$
It follows from
$$
     r_{\alpha \beta } + r_{\beta \alpha } + 1 | order_{\alpha  \beta } (g)
$$
and
$$
      order_{\alpha  \beta } (g) |  r_{\alpha \beta } + r_{\beta \alpha } + 1
$$
that
$$
          order_{\alpha  \beta } (g) = r_{\alpha \beta } + r_{\beta \alpha } + 1.
$$
Q.E.D.\\

Situations (i) and (ii)  from Theorem  \ref{morethan2} are illustrated by Fig.2 and Fig.3, respectively.

\begin{figure}[tb]
  \begin{center}
     \psfig{file=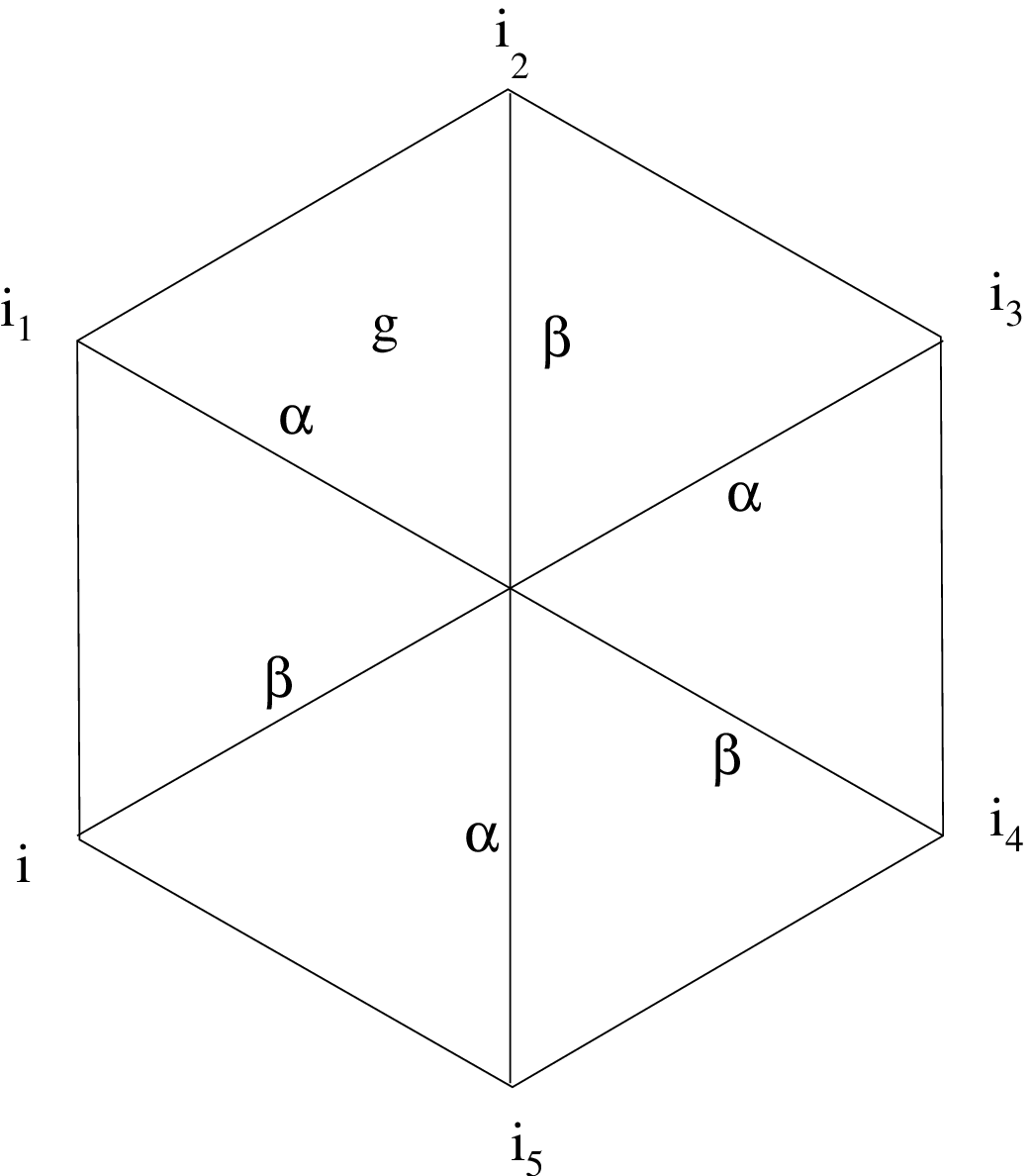,scale=0.7}         
      \caption{$(p_\alpha \cdot p_\beta )^3 (g) =g;$ the edges marked with $\alpha$ and $\beta$ are representatives of the corresponding equivalence classes.}
    \label{fig:order6} 
  \end{center}
\end{figure}

\begin{figure}[tb]
  \begin{center}
     \psfig{file=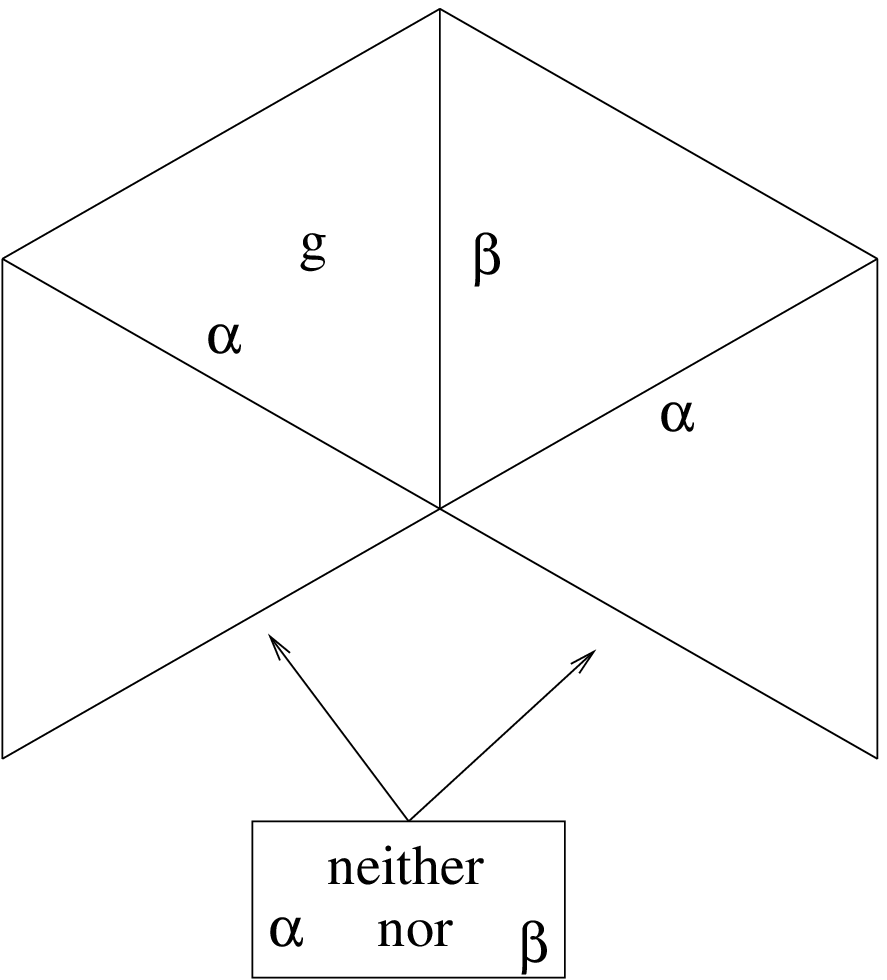,scale=0.7}         
      \caption{$(p_\alpha \cdot p_\beta )^4 (g) =g;$ the edges marked with $\alpha$ and $\beta$ are representatives of the corresponding equivalence classes.}
    \label{fig:order6} 
  \end{center}
\end{figure}

\begin{theorem}
\label{sing}
Let $M=a\star S/\simeq$ be a closed stellar $(n+1)$-dimensional manifold. For any two $(n-1)$-simplexes $\alpha, \; \beta \in S/\simeq$ the natural number $m_{\alpha \beta}$ denotes the least common multiple of natural numbers $\{order_{\alpha \beta }(g)\;;\; g \in S\},$ 
$$
(p_\alpha \cdot p_\beta )^{m_{\alpha \beta } } = 1 \mbox{ and } (p_\alpha \cdot p_\beta )^r \not= 1 \mbox{ for } 0 <r < m_{\alpha \beta } .
$$
\end{theorem}

This theorem is a direct corollary of Theorems \ref{order2} and \ref{morethan2}.\\
The geometrical properties of a stellar $(n+1)$-manifold $M=a\star (S/\simeq )$ find their representations among algebraic characteristics of $Cox(S/\simeq ).$  For example, the next statement characterizes the existence of collapsible $n$-simplexes in $S/\simeq .$

\begin{theorem}
\label{Alternation}
Let $M=a\star (S/\simeq )$ be a closed stellar $(n +1)$-manifold such that for any $\alpha \in S/\simeq $ one can find an even number $2\cdot m_\alpha$ so that
$$
(p_0 \cdot p_\alpha )^{2\cdot m_\alpha }  = 1 \mbox{ and } (p_0 \cdot p_\alpha )^r \not= 1 \mbox{ for } 0 <r < 2\cdot m_\alpha  .
$$
Then $S/\simeq $ does not contain any collapsible $n$-simplex if, and only if,
the group $<p_\alpha >_{\alpha \in S/\simeq } $  generated by $\{p_\alpha \}_{\alpha \in S/\simeq } $ is a subgroup of the alternating group $A_N,$ where  $N$ is the number of $n$-simplexes in $S.$
\end{theorem}
{\bf Proof}
If an $n$-simplex $h \in S/\simeq $ is collapsible then the corresponding equivalence class contains two $n$-simplexes $ g,\; p_\alpha (g) \in S$ where $\alpha $ is an $(n - 1)$ -simplex of both $g$ and $p_\alpha (g).$ Hence, $p_\alpha $ is a transposition and $<p_\alpha >_{\alpha \in S/\simeq }$ can not be a subgroup of the alternating group $A_N.$\\
On the other hand, if $S/\simeq $ does not contain any collapsible $n$-simplexes then, by Theorems \ref{order2}, \ref{orderM} each $\alpha $  contains $2\cdot m_\alpha $ $(n - 1)$-simplexes. Thus, $p_\alpha $ is an even permutation for any $\alpha \in S/\simeq $  and consequently $p_\alpha \in A_N.$ Thus, $<p_\alpha >_{\alpha \in S/\simeq } $ is  a subgroup of the alternating group $A_N.$
Q.E.D.\\

\section{Flat stellar manifolds}
\begin{figure}[tb]
  \begin{center}
     \psfig{file=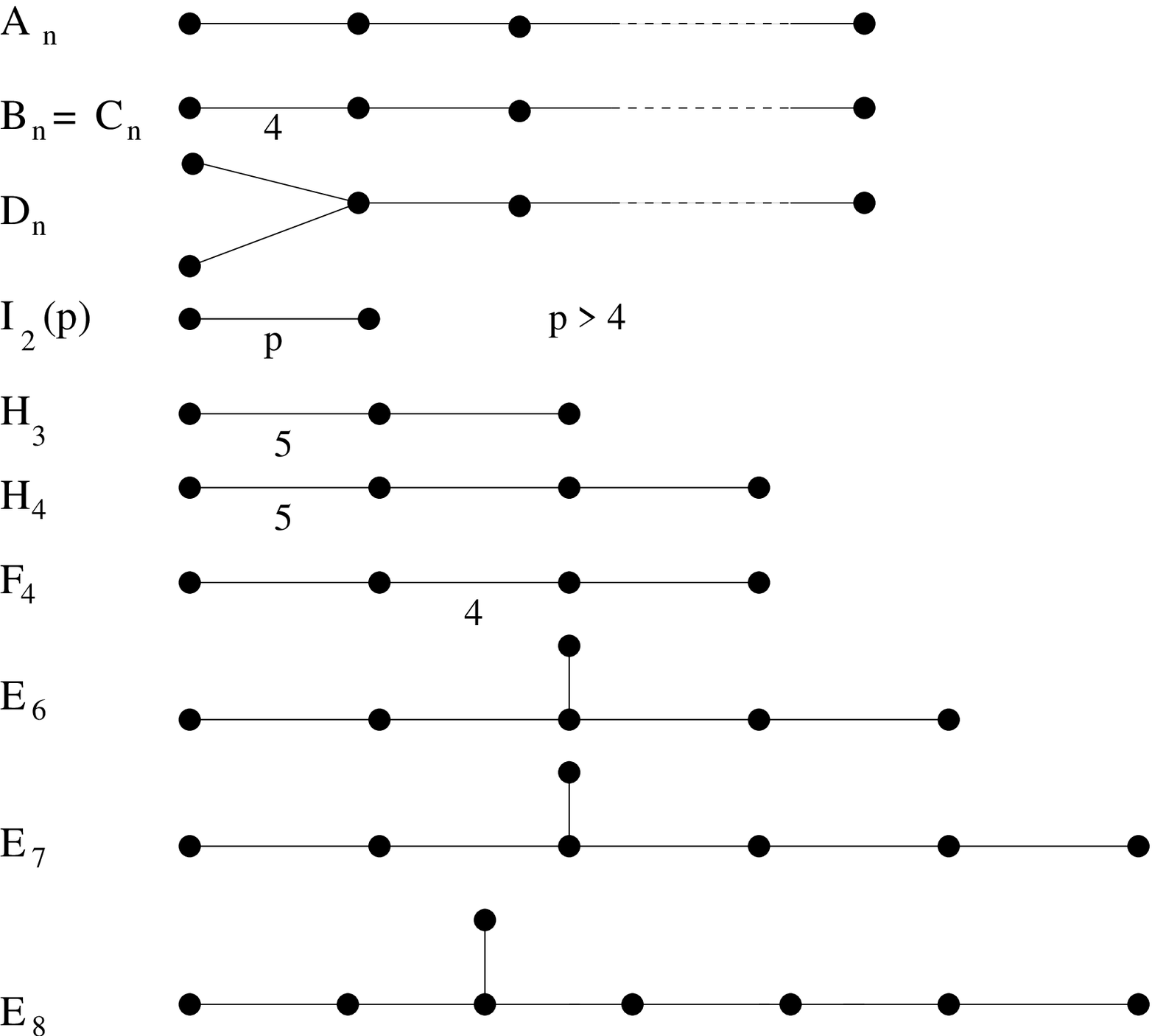,scale=0.7}         
      \caption{The classification of finite Coxeter groups}
    \label{fig:coxeter} 
  \end{center}
\end{figure}

Each stellar manifold $M$ has a stellar structure $M=a\star S/\simeq .$ Hence, each stellar manifold $M$ corresponds to a finite  Coxeter group $Cox(S/\simeq).$
Finite Coxeter groups are well studied and their classification is given in Fig.\ref{fig:coxeter}.

Each vertex  in Fig.\ref{fig:coxeter} represents a generator of a Coxeter group. If $m_{\alpha \beta } \ge 3$ then two vertices $\alpha$ and $\beta$ are connected by an arc. The arc is labeled with $m_{\alpha \beta }$ when $m_{\alpha \beta } > 3.$  Each Coxeter group is a direct product of the Coxeter groups corresponding to the connected graphs in Fig.\ref{fig:coxeter} and its diagram is the union of the diagrams of its factors. Thus,
$$
   Cox(S/\simeq ) = G_1  \otimes G_2 \otimes \dots \otimes G_N,
$$
where each $G_i \;(i=1,\;2,\dots,\; N)$ is one of the basic Coxeter groups corresponding to one of the diagrams in Fig.\ref{fig:coxeter}. \\
However, given a stellar structure $S/\simeq $ it might be not a trivial exercise to find out which set of Coxeter graphs in Fig. \ref{fig:coxeter} corresponds to $Cox(S/\simeq).$  The problem is that generators from  $\{p_0,\;\{p_\alpha \}_{\alpha \in S/\simeq } \}$ are not necessary independent. Thus one need to choose a system $S \subset \{p_0,\;\{p_\alpha \}_{\alpha \in S/\simeq } \} $ that consists of independent generators and  $Cox(S/\simeq)$ can be built out of $S.$ It is illustrated by the following example.\\

\begin{figure}[tb]
  \begin{center}
     \psfig{file=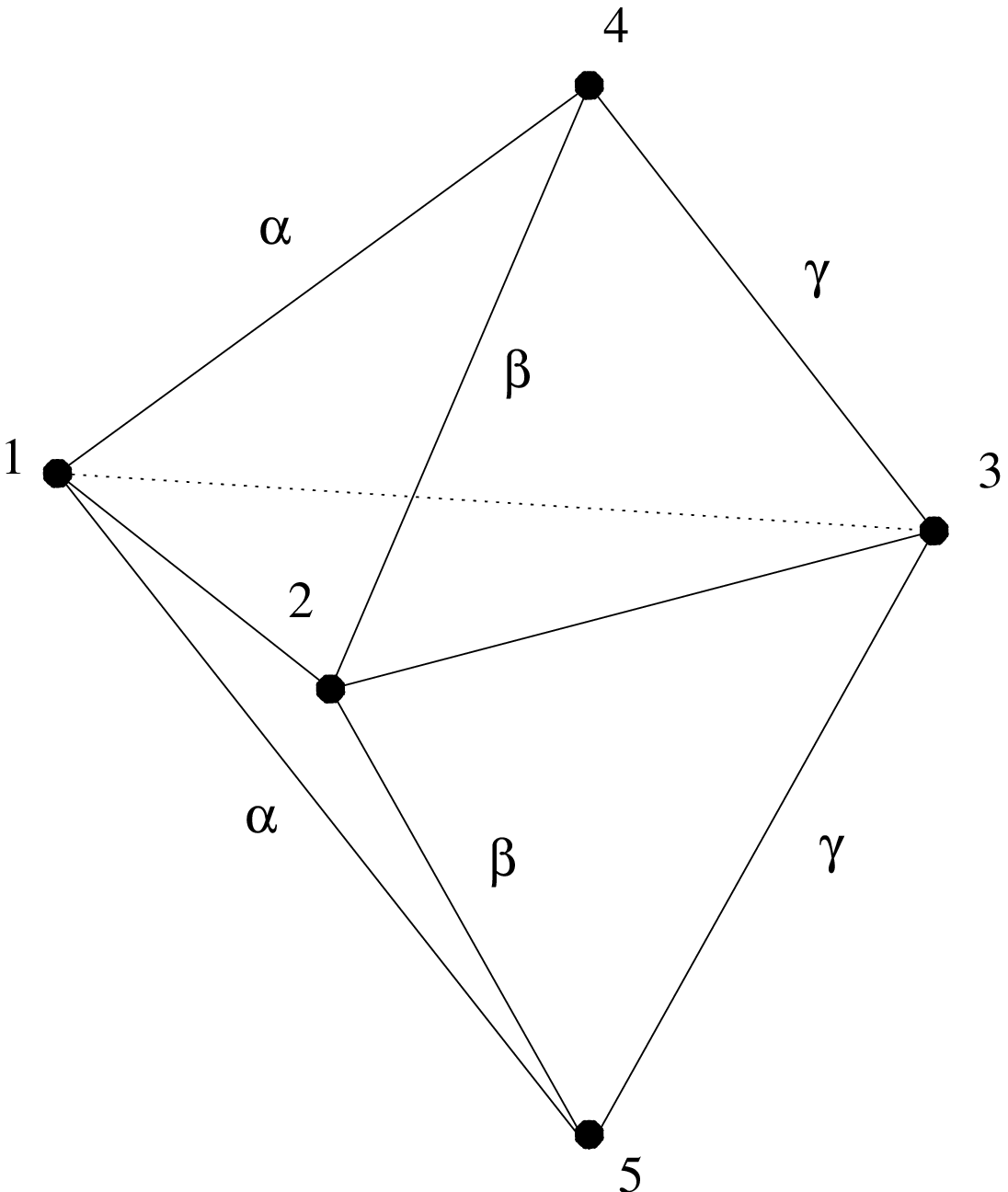,scale=0.7}         
      \caption{$S$}
    \label{sphere} 
  \end{center}
\end{figure}

\begin{figure}[tb]
  \begin{center}
     \psfig{file=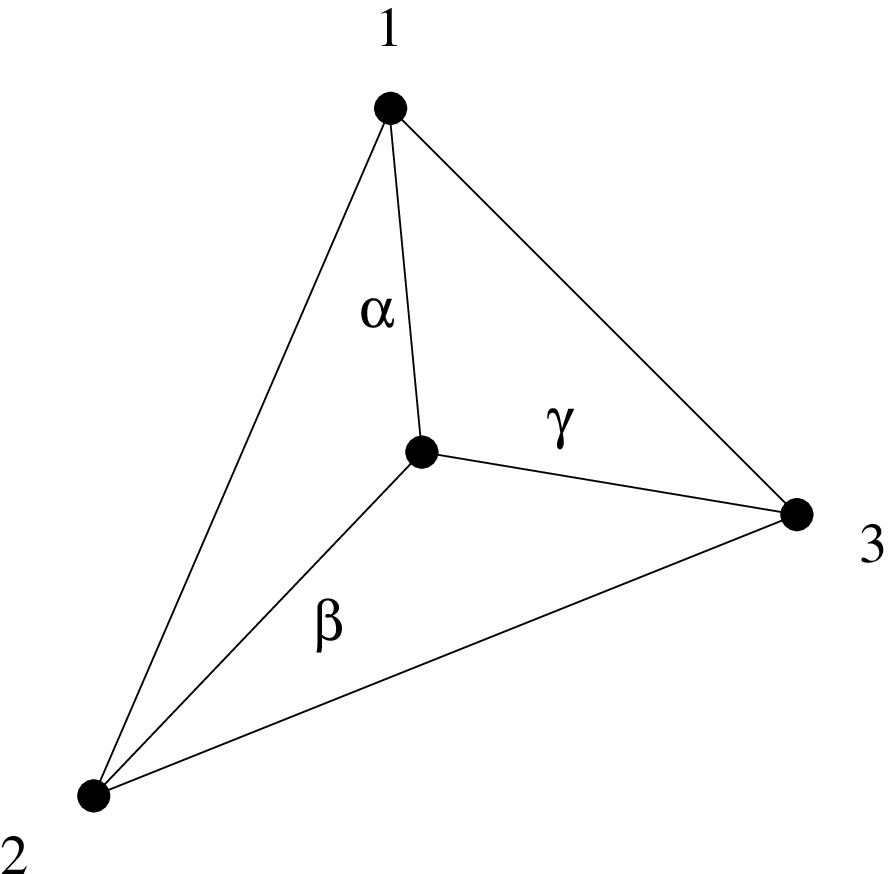,scale=0.7}         
      \caption{$S/\simeq $}
    \label{stellar} 
  \end{center}
\end{figure}

\begin{figure}[tb]
  \begin{center}
     \psfig{file=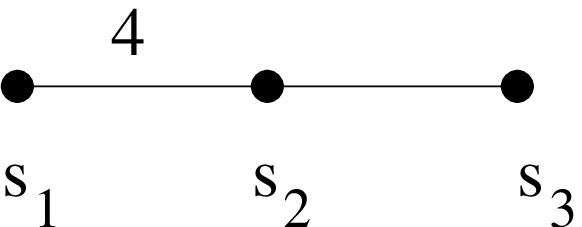,scale=0.7}         
      \caption{Diagram for $Cox(S/\simeq ).$}
    \label{diagram} 
  \end{center}
\end{figure}

{\bf Example 1.}
Consider a $2$-dimensional stellar sphere (Fig. \ref{sphere})
$$
S = (1\; 2 \; 4) + (1 \; 3 \; 4) + (2 \; 3 \; 4) + (1 \; 2 \; 5) + (1 \; 3 \; 5) + (2\; 3 \; 5).
$$
If regular equivalence is defined by $4 \simeq 5$ then the corresponding stellar structure is depicted in Fig. \ref{stellar}. It can be shown that Coxeter group $Cox(S/\simeq )$ has the following system of generators
$$
s_1 =  p_{( 1 \; 2) }, \;\; s_2 =  p_{\beta },\;\;s_3 = p_{\gamma},
$$
where $\beta $ is the equivalence class $\{ (2 \; 4),\; (2 \; 5) \}$  and $\gamma = \{ (3 \; 4),\; (3 \; 5) \} .$ 
The corresponding Coxeter diagram is depicted in Fig. \ref{diagram}.  Q.E.D.\\

All $1$-simplexes of the stellar structure from Example 1 are flat. Such stellar manifolds are called flat. 
\begin{definition}
A stellar manifold $M=a\star (S/\simeq) $ is called flat if, and only if, the group 
$$
\{1,\; p_0 \}
$$
is a normal subgroup of $Cox(S/\simeq).$
\end{definition} 

\begin{figure}[tb]
  \begin{center}
     \psfig{file=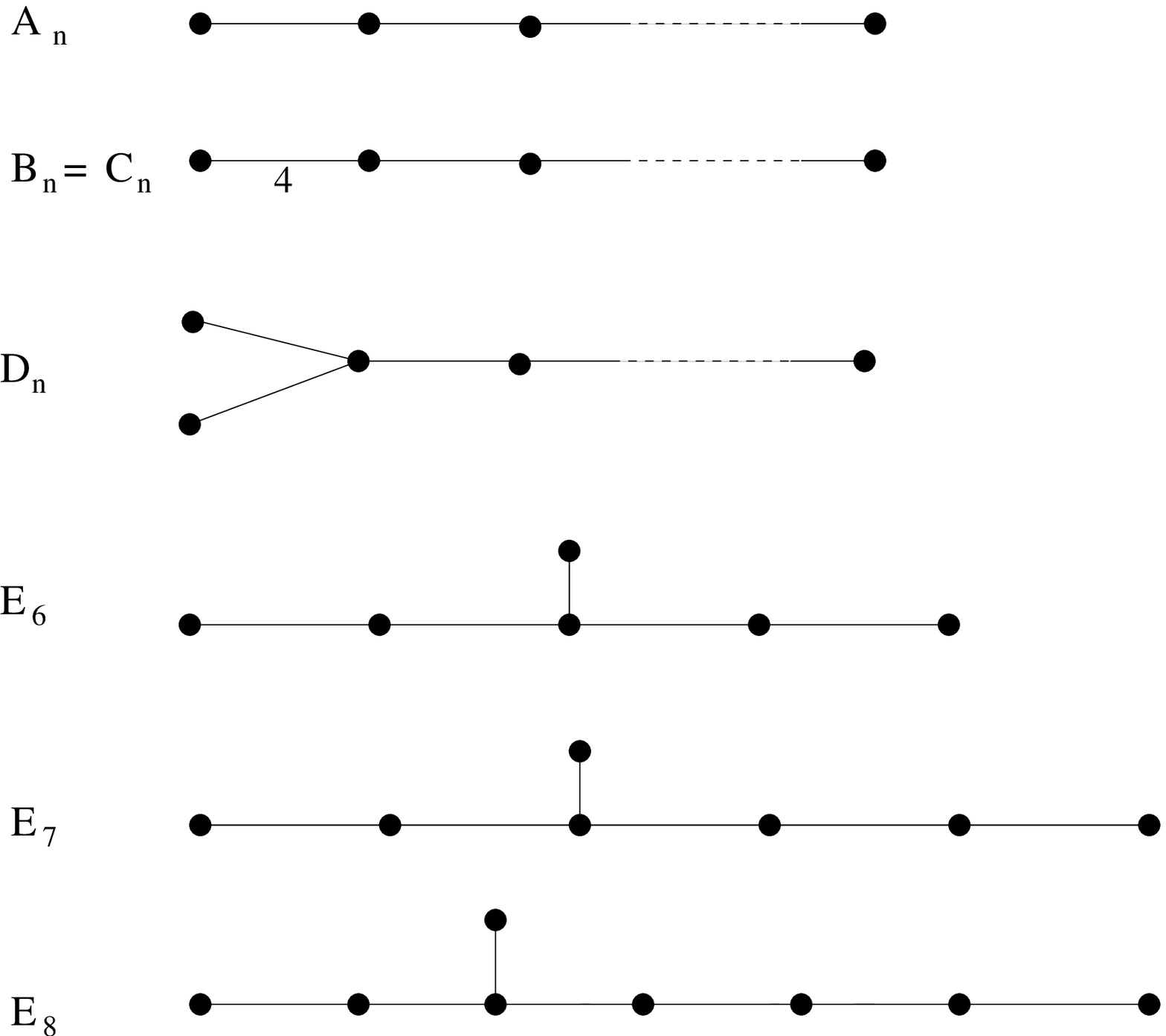,scale=0.7}         
      \caption{Coxeter diagrams corresponding to flat stellar manifolds}
    \label{flatCoxeter}
  \end{center}
\end{figure}

Our interest in flat stellar manifolds is motivated by the believe that an arbitrary stellar manifold can be built out of flat pieces glued together along singular simplexes defined by Theorem \ref{orderM}. Moreover, it follows from Theorems  \ref{order2}, \ref{orderM}, \ref{sing} that Coxeter group of any flat stellar manifold is a direct product of Coxeter groups with diagrams depicted in Fig. \ref{flatCoxeter} .

That means any flat stellar manifold can be build out of simple flat stellar manifolds corresponding to irreducible Coxeter diagrams in Fig. \ref{flatCoxeter}.\\
The next statement present a complete classification of flat $3$-dimensional stellar manifolds. 

\begin{theorem}({\bf Classification of flat $3$-manifolds})
\label{involution}
If $M$ is a closed connected flat stellar $3$-manifold then it is stellar equivalent to
$$
a\star (S/\simeq ), 
$$
where 
$$
S/\simeq = \left\{ \begin{array}{cc}
                         \mbox{disk } D_2 & \mbox{ for } \pi (M) = \{1\} \\
                         \mbox{projective plane } P_2 & \mbox{ for } \pi (M) = {\bf Z}_2
                  \end{array}
           \right.
$$
\end{theorem}
{\bf Proof.}
If $M$ is a flat stellar manifold  then $M=a\star (S/\simeq )$ and $ \{1,\; p_0 \} $ is a normal subgroup of $Cox(S/\simeq ).$
The involution $p_0$ corresponds to a piecewise linear one-to-one involution on $S^2.$ It follows from  \cite{Lopez} that there exists a circle $S^1$ that splits $S^2$ into two parts $S_+$ and $S_-$ such that
$$
p_0:\; S^1 \to S^1
$$
and $p_0$ is a one-to-one mapping that interchanges $S_+$ and $S_-.$
By Theorem \ref{Euler}
$$
\chi (S/\simeq) =1.
$$
Hence, $S/\simeq $ is either a disk $D_2$ or a projective plane $P_2$ (see, e.g., \cite{Massey} for details).
Q.E.D.\\
 Flat stellar manifolds have a range of properties that allow to tackle certain difficult problems. In particular, the next statement proves the famous Poincare conjecture for flat stellar manifolds.
 
\begin{theorem}({\bf Poincare conjecture for flat stellar manifolds}) 
If $M$ is compact, closed, connected and simply connected $3$-dimensional stellar flat manifold then $M$ is a stellar $3$-sphere. 
\end{theorem}   
{\bf Proof.}
If $M$ is a closed connected flat $3$-manifold and $\pi(M)=\{1\}$ then by Theorem \ref{involution}  $S/\simeq$ is a disk.
Since  $S/\simeq $ is collapsible then so is the prism $P(S/\simeq ).$  Due to \cite{Whitehead}  the prism $P(S/\simeq)$ is a stellar ball.
$M$ is stellar equivalent to 
$$
L=\prod_{b \in S/\simeq} ((a\; b)\; c_b) (a\star (S/\simeq )),
$$
where $\prod_{b \in S/\simeq} ((a\; b)\; c_b)$ is a superposition of subdivisions $((a\; b)\; c_b),$ where $c_b \notin S/\simeq$ and $c_b \not= c_d $ for $b \not= d.$
$$
 L=  a \star \partial P(S/\simeq ) + P(S/\simeq ),
$$
$L$ is a union of two balls with identified boundaries. Hence \cite{Glaser}, $L$ is a stellar $3$-sphere and so is $M.$
Q.E.D.\\

This simple proof is not valid for stellar manifolds that are not flat. The obstacles are created by singular simplexes defined as follows. \\

\begin{definition} ({\bf Singularity})
Let $M=a\star (S/\simeq )$ be a closed stellar $(n+1)$-dimensional manifold. An $(n-1)$-simplex $\alpha \in S/\simeq $ is called singular if the equivalence class $\alpha $ contains more than two different $(n-1)$-simplexes from $S.$
\end{definition}

By Theorem \ref{orderM} an $(n-1)$-simplex $\alpha \in  S/\simeq $ is singular only when there exists $m\ge 3$ such that
$$
(p_0 \cdot p_\alpha )^m = 1  \mbox{ and }  (p_0 \cdot p_\alpha )^{k} \not= 1 \mbox{ for } 0<k < m.
$$
Classification of Coxeter groups together with Theorems \ref{order2}, \ref{orderM}, \ref{sing} suggest the list of all possible singularities for a stellar structure. However, this topic is beyond the scope of this publication.

We conclude this paper by example showing that $\pi (M) = 1$ is essential for a flat $3$-manifold to be a $3$-sphere.\\

\begin{figure}[tb]
  \begin{center}
     \psfig{file=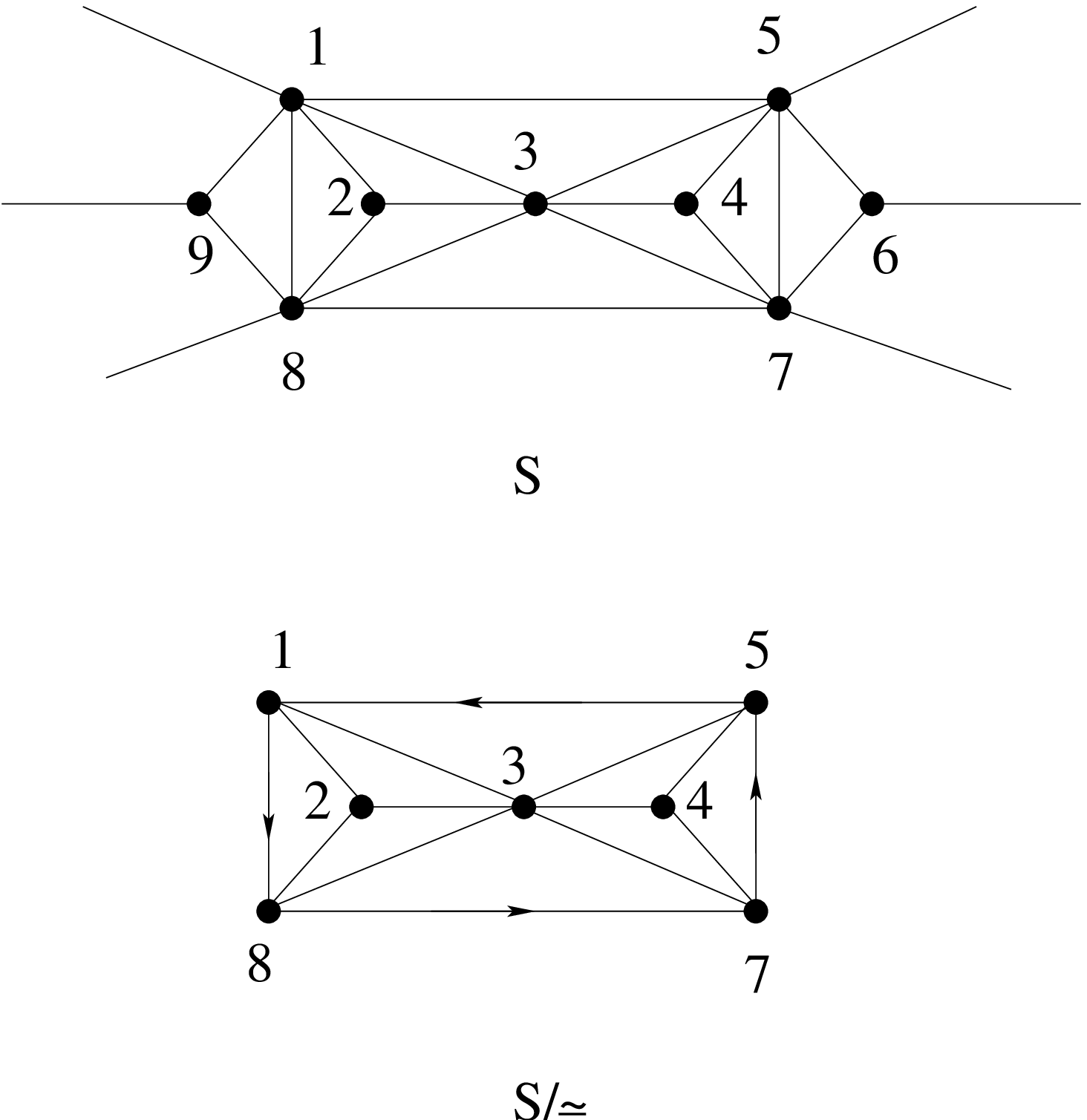,scale=0.7}         
      \caption{$S/\simeq $ is a projective plane.}
    \label{projectivePlane}
  \end{center}
\end{figure}

{\bf Example 2.}
Consider a triangulation of $2$-sphere $S$ depicted in Fig.\ref{projectivePlane}. Let us introduce regular equivalence so that
$$
2 \simeq 6,\;\;4 \simeq 9, \;\; 3 \simeq \infty ,\;\; 5\simeq 8, \;\; 1 \simeq 7
$$ 
where $\infty$ denotes the vertex at infinity (the south pole of $S$ while "$3$" is its north pole), and only the following generators are equivalent to each other:
$$
(3\; 7\; 8) \simeq (1\; 5 \; \infty ),\;\; (1\; 3 \; 5) \simeq (7 \; 8 \; \infty), \;\; (3\; 4 \; 5) \simeq (8 \; 9 \; \infty ),\;\; (3\; 4 \; 7) \simeq (1 \; 9 \; \infty ) ,
$$
$$
(4 \; 5 \; 7) \simeq (1\; 8 \; 9),\;\; (5\; 6 \; 7) \simeq (1\; 2 \; 8), \;\; (1\; 2 \; 3) \simeq (6 \; 7 \; \infty ),\;\; (2\; 3 \; 8) \simeq (5 \; 6 \; \infty ).
$$
 Then $S/\simeq $ is depicted in Fig.\ref{projectivePlane} where we need to glue the boundary so that $1\simeq 7$ and $5 \simeq 8 .$ Clearly, $ S/\simeq $ is a projective plane and $\chi (S/\simeq ) =1.$ It is easy to see  that $M=a\star (S/\simeq)$ is a stellar flat $3$-manifold. However, $\pi (M) = {\bf Z}_2 .$ 

Q.E.D.\\

\bibliographystyle{ams-plain}

\end{document}